\documentclass{mrlart2e}
\def\volno{12}
\def\yrno{2005}
\def\lpageno{100NN} 
\usepackage{latexsym}

\usepackage{amsmath,amssymb,amsfonts,amscd,amsthm}

\def\CC{\mathbb{C}}

\def\QQ{\mathbb{Q}}
\def\RR{\mathbb{R}}
\def\ZZ{\mathbb{Z}}

\def\og{\overline{\gamma}}
\def\ug{\underline{\gamma}}
\def\delbar{\overline{\partial}}

\author{Fr\'ed\'eric Bourgeois}
\title[Contact homology and homotopy groups]{Contact homology and homotopy groups \\
of the space of contact structures}
\date{}
\address{Universit\'e Libre de Bruxelles, D\'epartement de Math\'ematiques CP 218, Boulevard du Triomphe, 1050  Bruxelles, Belgium} 
\email{fbourgeo@ulb.ac.be}

\newtheorem{proposition}{Proposition}
\newtheorem{lemma}{Lemma}

\newtheorem*{conjecture}{Conjecture}

\theoremstyle{definition} 
\newtheorem{remark}{Remark}

\begin{document}

\begin{abstract}
Using contact homology, we reobtain some recent results of Geiges and Gonzalo about the 
fundamental group of the space of contact structures on some 3-manifolds. We show 
that our techniques can be used to study higher-dimensional contact manifolds and higher 
order homotopy groups.
\end{abstract}

\maketitle

\section{Introduction}

Geiges and Gonzalo \cite{GG} recently studied the topology of the space $\Xi(M)$ of 
contact structures for some $3$-manifolds $M$. In particular, they showed that the 
fundamental group $\pi_1(\Xi(T^3), \xi_n)$ based at any tight contact structure $\xi_n$ 
on $T^3$ contains an infinite cyclic group. They obtain a similar result in the more 
general case of $T^2$-bundles over $S^1$. Their argument involves classical techniques of
contact geometry, but relies on the classification of tight contact structures on these manifolds. 
It is therefore limited to $3$-dimensional contact manifolds.

The aim of this paper is to present an alternative proof for this result  that is based on 
contact homology. The advantage of this technique is that it can be used to study the 
fundamental group $\pi_1(\Xi(M),\xi_0)$ for higher-dimensional 
contact manifolds $M$, as 
well as homotopy groups of higher order $\pi_k(\Xi(M),\xi_0)$.

Contact homology is an invariant of contact structures that was recently introduced by 
Eliashberg, Givental and Hofer \cite{EGH}. The contact homology $HC_*(M,\xi)$ of a 
contact structure $\xi$ on $M$ is defined as the homology of a chain complex 
$(C_*,d)$. Let $\alpha$ be a contact form for $\xi$, i.e. $\xi = \ker \alpha$. 
The Reeb vector field $R_\alpha$ associated to $\alpha$ is characterized by 
$\imath(R_\alpha) d\alpha = 0$ and $\alpha(R_\alpha) = 1$. It is possible to choose 
$\alpha$ generically, so that the closed orbits of $R_\alpha$ are nondegenerate; in 
particular, the set of closed Reeb orbits is countable. For simplicity, we assume that 
the Reeb vector field has no contractible closed orbits. Fix a nontrivial free homotopy class ${\bar a}$
in $M$ and let $C^{\bar a}_*$ be the module that is freely generated by the closed Reeb 
orbits in the free homotopy class ${\bar a}$, over the coefficient
ring $\QQ[H_2(M,\ZZ)/\mathcal{R}]$, where $\mathcal{R}$ is a submodule of $H_2(M,\ZZ)$
such that $\mathcal{R} \subset \ker c_1(\xi)$. 

We choose a compatible complex structure on $\xi$ and extend it 
to an almost complex structure $J$ on the symplectization $(\RR \times M, d(e^t \alpha))$
by $J \frac{\partial}{\partial t} = R_\alpha$. We consider those $J$-holomorphic cylinders in
this symplectic manifold that converge to closed Reeb orbits as $t \to \pm \infty$. We can
associate an element of $H_2(M,\ZZ)/\mathcal{R}$ to such a $J$-holomorphic cylinder
in the following way. The homotopy class ${\bar a}$ being nontrivial, we choose a loop  
representing ${\bar a}$ and a homotopy between that loop and any closed Reeb orbit in
class ${\bar a}$. Attaching the projection to $M$ of a $J$-holomorphic cylinder to a pair of 
the above homotopies, we obtain a $2$-torus in $M$ and hence a homology class in
$H_2(M,\ZZ)$. 

We then define a linear morphism $d : C^{\bar a}_* \to C^{\bar a}_{*-1}$ of degree $-1$ that 
counts the $J$-holomorphic cylinders of index $1$ in $(\RR \times M, d(e^t \alpha))$ and that encodes 
their homology class using the coefficient ring. The morphism $d$ has the property that $d \circ d = 0$
so that $(C_*^{\bar a}, d)$ is a chain complex. This follows from the fact that $d \circ d$ counts 
elements in the boundary of the one-dimensional moduli space of $J$-holomorphic cylinders of index 
$2$. After orienting the moduli spaces \cite{BM} and associating a sign \cite{FH} to each 
$J$-holomorphic cylinder of index $1$, two terms in $d \circ d$ corresponding to boundary points 
of the same connected component in the moduli space will cancel each other. The contact homology 
$HC_*^{\bar a}(M,\xi)$ of $(M,\xi)$ is then defined as the homology of the chain complex 
$(C_*^{\bar a}, d)$. It is independent of the choice of a contact form $\alpha$ and of a 
complex structure $J$ for $\xi$. 

We will show how to use contact homology to study the space of contact structures.
Following the philosophy of this paper, K\'alm\'an \cite{K} obtained similar results for the fundamental
group of the space of Legendrian embeddings in dimension $3$, using relative contact homology. 
\\[0.3cm]
{\it Acknowledgements. } This work was initiated during an EDGE postdoc at Ecole Polytechnique.
The author is grateful to Emmanuel Giroux and the referee for suggesting improvements to this paper.

\section{Automorphisms of contact homology}  \label{sec:pi1}

In view of Gray stability, the primary use of contact homology is to study the set 
$\pi_0(\Xi(M))$ of connected components of the space of contact structures. In this 
section, we explain how to study $\pi_1(\Xi(M),\xi_0)$ using automorphisms of contact homology.

Let $\xi_s$, $s \in \RR/\ZZ$, be a loop of contact structures based at $\xi_0$. Choose a
loop of contact forms $\alpha_s$  such that $\xi_s = \ker \alpha_s$ for 
$s \in \mathbb{R}/\mathbb{Z}$.  Let $f : \RR \to [0,1]$ be a smooth 
decreasing function such that $f(t) = 1$ if $t$ is sufficiently small and $f(t) = 0$ 
if $t$ is sufficiently large. If the function $f$ varies slowly enough, 
the $2$-form $d(e^t \alpha_{f(t)})$ on $\RR \times M$ is symplectic; 
we denote the corresponding symplectic cobordism by $(\RR \times M, \omega)$.

Let $J$ be a compatible almost complex structure on $(\RR \times M, \omega)$ such that,
for $|t|$ sufficiently large, $J \frac{\partial}{\partial t} = R_{\alpha_0}$, 
$J\xi_0 = \xi_0$ and $J$ is independent of $|t|$. The choice of the complex structure $J$ on $\xi$
for $|t|$ large is dictated by the following Lemma.
 
\begin{lemma} \label{lem:transversal}
Let $(M,\xi)$ be a contact manifold with contact form $\alpha$ and assume that the Reeb field
$R_\alpha$ has no contractible closed orbits. Let ${\bar a}$ be a free homotopy class containing 
only simple closed Reeb orbits.

Then, for a generic $\RR$-invariant compactible almost complex structure $J$ on the symplectization 
$(\RR \times M, d(e^t\alpha))$ satisfying $J\xi = \xi$ and $J\frac{\partial}{\partial t} = R_\alpha$,
all $J$-holomorphic cylinders asymptotic to closed Reeb orbits in class ${\bar a}$ are generic, i.e. 
their linearized Cauchy-Riemann operator is surjective.

\end{lemma}

In other words, under the assumptions of Lemma \ref{lem:transversal},
transversality can be achieved rather easily. Note that a similar result clearly holds in the
case of symplectic cobordisms, because the almost complex structure $J$ does not
satisfy any other condition than compatibility with the symplectic structure away from the ends.

Lemma \ref{lem:transversal} was originally proved by Dragnev \cite{D}, and will be proved 
independently in the Appendix.

Let $\Phi : C_* \to C_*$ be the linear morphism obtained by counting $J$-holomorphic 
cylinders of index $0$ in $(\RR \times M, \omega)$.

\begin{proposition} \label{prop:aut}
The above construction induces a group morphism
$$
\eta : \pi_1(\Xi(M), \xi_0) \to {\rm Aut}(HC_*(M,\xi_0)) .
$$
\end{proposition}

\begin{proof}
The linear morphism $\Phi$ is a chain map, i.e. satisfies the relation 
$\Phi \circ d = d \circ \Phi$. This follows from the fact that $\Phi \circ d - d \circ \Phi$
counts elements in the boundary of the one-dimensional moduli space of $J$-holomorphic cylinders 
of index $1$ in $(\RR \times M, \omega)$. As in the case of $d \circ d$, all terms will cancel 
in pairs. Therefore, $\Phi$ induces an endomorphism 
$\overline{\Phi} : HC_*(M,\xi_0) \to HC_*(M,\xi_0)$. In order to show that 
$\overline{\Phi}$ is an isomorphism, we construct an inverse of $\overline{\Phi}$. 
Consider the reversed loop $\xi'_s = \xi_{1-s}$, for $s \in [0,1]$. We write 
$(\RR \times M,\omega')$ for the symplectic cobordism constructed from $\xi'$ as above, 
and $\overline{\Phi}'$ for the morphism of contact homology corresponding to this 
symplectic cobordism. Note that we can glue the symplectic cobordisms 
$(\RR \times M, \omega)$ and $(\RR \times M, \omega')$ and obtain a symplectic cobordism 
$(\RR \times M, \omega \sharp \omega')$. The morphism of contact homology corresponding
to this symplectic cobordism is the composition $\overline{\Phi}' \circ \overline{\Phi}$.
The symplectic form $\omega \sharp \omega'$ can be deformed into the symplectic form 
$d(e^t \alpha_0)$. Since the contact homology morphism corresponding to a symplectization
is the identity, it follows that $\overline{\Phi}'$ is a left inverse for
$\overline{\Phi}$. We prove similarly that $\overline{\Phi}'$ is a right inverse for 
$\overline{\Phi}$. 

We set $\eta(\xi_s) = \overline{\Phi}$. We need to check that $\eta$ is well-defined, 
i.e. that $\overline{\Phi}$ is invariant under deformation of $\xi_s$. Given a 
$2$-parameter family $\xi_{s,t}$, for $s,t \in [0,1]$ and satisfying 
$\xi_{0,t} = \xi_{1,t} = \xi_0$ for $t \in [0,1]$, we can construct a $1$-parameter 
family of symplectic cobordisms $(\RR \times M, \omega_t)$ equipped with compatible 
almost complex structures $J_t$. The chain maps $\Phi_i : (C_*,d) \to (C_*,d)$, 
$i = 0,1$ corresponding to the symplectic cobordisms at $t = 0$ and $t=1$ satisfy
the relation $\Phi_1 - \Phi_0 = K \circ d + d \circ K$, where $K : C_* \to C_{*+1}$ is a 
linear morphism of degree $1$ counting $J_t$-holomorphic cylinders of degree $-1$ in 
$(\RR \times M, \omega_t)$, for $t \in [0,1]$. This relation is obtained by counting 
elements in the boundary of the one-dimensional moduli space of $J_t$-holomorphic
cylinders of index $0$ in $(\RR \times M, \omega_t)$ for $0 \le t \le 1$. 
It follows that the induced isomorphisms $\overline{\Phi}_0$ and $\overline{\Phi}_1$ 
on contact homology are identical.

Finally, we show that the map $\eta$ is a group morphism. Given $1$-parameter families 
$\xi_s$ and $\xi'_s$, for $s \in [0,1]$, such that $\xi_0 = \xi_1 = \xi'_0 = \xi'_1$, 
we can construct symplectic cobordisms $(\RR \times M, \omega)$ and 
$(\RR \times M, \omega')$. Then, the glued cobordism 
$(\RR \times M, \omega \sharp \omega')$ corresponds, up to deformation, to the 
catenation of the paths $\xi_s$ and $\xi'_s$. But we already saw that the 
isomorphism $\eta([\xi_s].[\xi'_s])$ counting $J$-holomorphic cylinders of index $0$ in 
$(\RR \times M, \omega \sharp \omega')$ is the composition $\eta(\xi_s) \circ \eta(\xi'_s)$
of the isomorphisms $\eta(\xi_s)$ and $\eta(\xi'_s)$.
\end{proof}

This proposition shows in particular that a loop $\xi_s$ in $\Xi(M)$ 
induces a nonzero element of $\pi_1(\Xi(M),\xi_0)$ as soon as $\eta(\xi_s) \neq id$.

We now turn to a particular case for which it is very easy to compute the automorphism 
$\eta(\xi_s)$.  Assume that there exists a loop $\phi_s$, $s \in \mathbb{R}/\mathbb{Z}$,
of diffeomorphisms of $M$ such that $\xi_s = \phi_{s*} \xi_0$ for all 
$s \in \mathbb{R}/\mathbb{Z}$. Let $\alpha_0$ be a contact form for $\xi_0$ and
$J_0$ be an $\RR$-invariant almost complex structure in the symplectization
$(\RR \times M, d(e^t \alpha_0))$. Then we may choose the symplectic cobordism
$(\RR \times M,\omega)$ in the definition of $\eta(\xi_s)$ to be the image of the
symplectization $(\RR \times M,d(e^t \alpha_0))$ under the diffeomorphism
$(t,p) \to (t,\phi_{f(t)}(p))$. Moreover, we may choose the almost complex structure $J$
on $(\RR \times M,\omega)$ to be the image of $J_0$ under the above diffeomorphism. 

Now, if $J_0$ is generic as in Lemma \ref{lem:transversal}, the only $J_0$-holomorphic 
cylinders in the symplectization are the vertical ones. Hence, the only $J$-holomorphic 
cylinders in $(\RR \times M,\omega)$ are the cylinders 
$(t,\vartheta) \to (t, \phi_{f(t)} \circ \gamma(\vartheta))$. We can summarize these
observations in the following lemma.

\begin{lemma} \label{lem:holcyl}
Let ${\bar a}$ be a free homotopy class containing only simple closed Reeb orbits. 
Suppose that $\xi_s = \phi_{s*} \xi_0$ for all 
$s \in \mathbb{R}/\mathbb{Z}$ where $\phi_s$, $s \in 
\mathbb{R}/\mathbb{Z}$, is a loop of diffeomorphisms of $M$.
Then, for any closed Reeb orbit $\gamma$ in class ${\bar a}$, 
$$
\eta(\xi_s) \gamma = e^A \gamma ,
$$
where $A \in H_2(M,\ZZ)$ is the homology class of the torus
$(t,\vartheta) \to \phi_{f(t)} \circ \gamma(\vartheta)$. 
\end{lemma}

Note that Lemma \ref{lem:holcyl} also holds if $\phi_s$, $s \in [0,1]$, is a path of diffeomorphisms
such that $\phi_1$ preserves the generic contact form $\alpha$ that is used to compute the contact
homology of $\xi_0$.

\section{Fundamental group of the space of contact structures}

In this section, we reprove the result of Geiges and Gonzalo \cite{GG} about 
the fundamental group of the space of contact structures on $T^2$-bundles over
$S^1$ and extend it to higher dimension.

We first consider the case of the $3$-torus. Recall (see for example \cite{Kanda}) that any tight contact structure $\xi$ on $T^3$ is diffeomorphic to one of the contact structures $\xi_n = \ker \alpha_n$, 
$n = 1, 2, \ldots$, where
$$
\alpha_n = \cos(n\theta) dx + \sin(n\theta) dy .
$$

\begin{proposition}  \label{prop:T3}
The fundamental group $\pi_1(\Xi(T^3),\xi_n)$ based at any tight contact structure $\xi_n$ 
on $T^3$ contains an infinite cyclic subgroup.
\end{proposition}

\begin{proof}
Consider the loop of contact structures $\xi_{n,s} = \ker \alpha_{n,s}$ where 
$$
\alpha_{n,s} = \cos(n\theta - 2\pi s)dx + \sin(n\theta - 2\pi s) dy .
$$

Consider the homotopy class ${\bar a}$ of loops wrapping once along $T^3$ in the 
$x$-direction. There are $n$ circles of closed Reeb orbits in the class ${\bar a}$, at angles
$\theta = \frac{2\pi}{n}k$, $k = 0, \ldots, n-1$.
Let $HC^{\bar a}_*(T^3,\xi_n)$ be the contact homology over the coefficient ring
$\QQ[H_2(T^3,\ZZ)/\mathcal{R}]$, where $\mathcal{R}$ is generated by the torus
$A_{x,y} \in H_2(T^3,\ZZ)$ along the $x$- and $y$-coordinates.
Each circle gives two generators, see \cite{Bproceed}, in $HC^{\bar a}_*(T^3,\xi_n)$ with 
degree $-1$ and $0$. Let $\gamma_k$ be the closed orbit of degree $0$ at 
$\theta = \frac{2\pi}{n}k$.

By Lemma \ref{lem:holcyl} and the remark below it, the image of $\gamma_k$ under the automorphism 
$\eta(\xi_{n,s})$ is proportional to $\gamma_{k+1}$.
Moreover,  if we choose $\gamma_0$ as representative of ${\bar a}$ and cylinders
$(t,\vartheta) \to (x_0 + \vartheta, y_0, \frac{2\pi}{n}k \, t)$ as homotopies from $\gamma_0$
to $\gamma_k$, $k = 0, \ldots, n-1$, the above $J$-holomorphic cylinder has a trivial
homology class except if $k = n-1$. In this last case, we obtain the homology class 
$A_{x,\theta} \in H_2(T^3,\ZZ)$ of the $2$-torus spanned by the direction of $x$ and $\theta$. 

Hence, on $HC_0^{\bar a}(T^3,\xi_n)$, the automorphism $\eta(\xi_{n,s})$ is given by
$$
\eta(\xi_{n,s}) \gamma_k = \left\{ \begin{array}{ll}
\gamma_{k+1} & \textrm{if } k \neq n-1 \\
\gamma_0 e^{A_{x,\theta}} & \textrm{if } k = n-1
\end{array} \right.
$$
Since all powers of this automorphism are distinct from the identity, the loop 
$\xi_{n,s}$ of contact structures generates an infinite cyclic subgroup in 
$\pi_1(\Xi(T^3),\xi_n)$.
\end{proof}

\begin{remark} The loop of contact structures induced by $\phi_s(x,y,\theta) = (x+s,y,\theta)$ or $\phi_s(x,y,\theta) = (x,y+s,\theta)$ is homotopically trivial, since the contact structures $\xi_n$ are $T^2$-invariant. Using our techniques, if we keep the coefficient ring $\QQ[H_2(T^3,\ZZ)/\mathcal{R}]$,
the contact homology automorphism induced by $\phi_s$ is the identity.
If instead we consider contact homology over the coefficient ring $\QQ[H_2(T^3,\ZZ)]$, then
$HC^{\bar a}_*(T^3,\xi_n)$ has $n$ generators $\gamma_1, \ldots, \gamma_n$ 
in degree $-1$, satisfying the relation $(1-e^{A_{x,y}}) \gamma_i = 0$. Therefore, multiplication
by $e^{A_{x,y}}$ is trivial and the contact homology automorphism induced by $\phi_s$ is 
the identity, as expected.
\end{remark}

We can extend this result to the more general case of $T^2$-bundles over $S^1$. 
For $A \in SL(2,\ZZ)$, let $T^3_A$ be the quotient of $T^2 \times \RR$ under the 
diffeomorphism
$$
\widetilde{A} : \Big( \left( \begin{array}{c} x \\ y \end{array} \right) , \theta \Big) \to 
\Big( A \left( \begin{array}{c} x \\ y \end{array} \right) , \theta + 2\pi \Big) .
$$
The manifold $T^3_A$ carries an infinite family $\zeta_n$, $n = 1, 2, \ldots$ of 
contact structures \cite{G} induced by $\widetilde{A}$-invariant contact structures on
$T^2 \times \RR$ of the form
$$
\cos f(\theta) dx + \sin f(\theta) dy = 0
$$
satisfying
$$
2(n-1)\pi < \sup_{\theta \in S^1} (f(\theta + 2\pi) - f(\theta)) \le 2n\pi .
$$
We can define loops of contact structures $\zeta_{n,s}$, $s \in \RR/\ZZ$, by the equation
$$
\cos f(\theta-2\pi s) dx + \sin f(\theta-2\pi s) dy = 0 .
$$

\begin{proposition}
The fundamental group $\pi_1(\Xi(T^3_A),\zeta_n)$ based at the contact structure $\zeta_n$ 
on $T^3_A$ contains an infinite cyclic subgroup.
\end{proposition}

\begin{proof}
With the above contact form for $\zeta_n$, each closed Reeb orbit is contained in a $T^2$-fiber
over  $\theta$ corresponding to a rational slope. Moreover, as in the case of $T^3$, each such 
fiber gives one generator of contact homology in degree $-1$ and one generator in degree $0$,
in the free homotopy class corresponding to $\theta$.

If $A$ has finite order in $SL(2,\ZZ)$, equivalently if $A$ is elliptic ($| \textrm{tr } A | < 2$) or 
$A = \pm I$, then any nontrivial, free homotopy class ${\bar a}$ in a $T^2$-fiber contains 
finitely many closed Reeb orbits. We are then in a situation analogous to Proposition 
\ref{prop:T3} and $\eta(\zeta_{n,s})$ is a cyclic permutation with multiplication of one generator 
by $e^{[T]} \neq 1$ in $H_2(T^3_A,\ZZ)$.

If $A$ is hyperbolic ($|\textrm{tr } A| > 2$), then any nontrivial, free homotopy 
class in a $T^2$-fiber contains infinitely many closed Reeb orbits $\gamma_k$, $k\in \ZZ$, 
of degree $0$. In this case, $\eta(\zeta_{n,s}) \gamma_k = \gamma_{k+1}$ for all $k \in \ZZ$. 

Finally, the case in which $A$ is parabolic ($|\textrm{tr } A| = 2$) but $A \neq \pm I$ is 
intermediate because the number of closed Reeb orbits in a free homotopy class 
${\bar a}$ may be finite or infinite, depending on ${\bar a}$. Therefore, the automorphism 
$\eta(\zeta_{n,s})$ is either a cyclic permutation with multiplication of one generator 
by a nontrivial unit or a shift,  depending on ${\bar a}$. 
\end{proof}

We now turn to the unit cotangent bundle $ST^*\Sigma_g$ of a compact, oriented surface 
$\Sigma_g$ of genus $g > 1$. Given a Riemannian metric on $\Sigma_g$, this manifold is 
naturally equipped with the Liouville $1$-form, defining a contact structure $\xi_{can}$.

\begin{proposition} \label{prop:STSigma}
The fundamental group $\pi_1(\Xi(ST^*\Sigma_g),\xi_{can})$ based at the contact 
structure $\xi_{can}$ on $ST^*\Sigma_g$ contains an infinite cyclic subgroup.
\end{proposition}

\begin{proof}
Let $\alpha_{can}$ be the Liouville $1$-form corresponding to a hyperbolic metric on $\Sigma_g$.
Since the Reeb flow on $ST^*\Sigma_g$ coincides with the geodesic flow, there is a unique
closed Reeb orbit in a fixed, nontrivial,  free homotopy class ${\bar a}$.

The principal $S^1$-bundle $ST^*\Sigma_g$ is naturally equipped with a fiberwise $S^1$-action
given by diffeomorphisms $\phi_s$, $s \in \RR/\ZZ$, of $ST^*\Sigma_g$. Let $\alpha_s = 
\phi_s^* \alpha_{can}$ and $\xi_s = \ker \alpha_s$. The corresponding loop of closed Reeb orbits
in class ${\bar a}$ spans the preimage in $ST^*\Sigma_g$ of the geodesic in $\Sigma_g$.
This is a torus $T$ that generates a free submodule in $H_2(ST^*\Sigma_g,\ZZ)$. Indeed, we can
find a loop in $\Sigma_g$ intersecting our closed geodesic exactly once; the natural lift of this
loop then has intersection number $1$ with $T$.

Hence, $\eta(\xi_s)$ is multiplication by $e^{[T]} \neq 1$, and $\xi_s$ generates an infinite
cyclic subgroup of  $\pi_1(\Xi(ST^*\Sigma_g),\xi_{can})$.
\end{proof}

We now extend these 3-dimensional results to the contact manifold $T^5$. Lutz \cite{L} gave 
the first example of a contact form on $T^5$,
using a knotted fibration $\varphi = (\varphi_1,\varphi_2) : T^3 \to \RR^2$ of the $3$-torus, 
given by
\begin{eqnarray*}
\varphi_1 &=& \epsilon (\sin \theta_1 \cos \theta_3 - \sin \theta_2 \sin \theta_3) ,\\
\varphi_2 &=& \epsilon (\sin \theta_1 \sin \theta_3 + \sin \theta_2 \cos \theta_3) .
\end{eqnarray*}

The corresponding submersion $\psi = \varphi / \| \varphi \| : T^3 \setminus \{ \varphi = 0 \} \to S^1$
defines an open book decomposition of $T^3$ with binding $\varphi^{-1}(0)$ and pages
$\psi^{-1}(p)$, where $p$ is a point in $S^1 \subset \RR^2$.
If $\epsilon > 0$ is sufficiently small, and if $\beta$ is a $1$-form on $T^3$ that does not 
vanish on the binding and such that $d\beta$ is an area form on each page, then the
$T^2$-invariant $1$-form
$$
\alpha = \varphi_1 d\theta_4 +
\varphi_2 d\theta_5 + \beta
$$
defines a contact structure $\xi = \ker \alpha$ on $T^5$. The Reeb field for
that contact form is proportional to
$$
\varphi_1 \frac{\partial}{\partial \theta_4} + \varphi_2 \frac{\partial}{\partial \theta_5}
+ X ,
$$
where $X$ is the Hamiltonian vector field of the function $\frac12(\varphi_1^2 + \varphi_2^2)$
restricted to each page, with respect to the symplectic form $d\beta$.
Let ${\bar a}$ be the free homotopy class of loops wrapping once along the $\theta_4$ coordinate.
The closed Reeb orbits in class ${\bar a}$ correspond to critical points of $\varphi_1^2$ on
the page $\psi^{-1}(1,0)$ on $T^3$. There are four maxima and eight saddles. The grading of the 
closed orbits for the maxima exceeds by $1$ the grading of the closed orbits for the saddles.

\begin{proposition}  \label{prop:T5}
The fundamental group $\pi_1(\Xi(T^5),\xi)$, based at the contact structure 
$\xi$, contains the subgroup $\ZZ^3$.
\end{proposition}

\begin{proof}
Let $\phi^{(i)}_s$, $i = 1, 2, 3$, be the loops of diffeomorphisms of $T^5$ defined by 
$\phi^{(i)}_s(\theta_1,\ldots ,\theta_5) = (\theta'_1,\ldots, \theta'_5)$ with 
$\theta'_i = \theta_i + s$ and $\theta'_j = \theta_j$ for $j \neq i$.
 Since the contact 
structure $\xi$ is invariant under $\phi^{(4)}_s$ and $\phi^{(5)}_s$, we can
consider $3$ circles of contact structures $\xi^{(i)}_s = (\phi^{(i)}_s)_* \xi$ for 
$i = 1, 2, 3$. We wish to compute $\eta(\xi^{(i)}_s)$ for $i = 1, 2, 3$.

Let $HC^{\bar a}_*(T^5,\xi)$ be the contact homology over the ring
$\QQ[H_2(T^5,\ZZ)/\mathcal{R}]$, where $\mathcal{R}$ is generated by 
the torus $A_{4,5} \in H_2(T^5,\ZZ)$ along the $\theta_4$- and $\theta_5$-coordinates. 
Consider a closed Reeb orbit $\gamma$ in class ${\bar a}$.

As in Proposition \ref{prop:T3}, the image of $\gamma$ under the automorphism
$\eta(\xi^{(i)}_s)$ is $\gamma e^{A_{i,4}}$, where $A_{i,4} \in H_2(T^5,\ZZ)$ is the 
homology class of the $2$-torus spanned by the directions of $\theta_i$ and $\theta_4$.

On the other hand, the closed Reeb orbits corresponding to the saddles induce in
$HC^{\bar a}_*(T^5,\xi)$ the quotient of a free module of rank $8$ by a free module
of rank at most $4$. Therefore, it contains at least a free summand. The automorphism
$\eta(\xi^{(i)}_s)$ on that summand is multiplication by $e^{A_{i,4}}$.

Therefore, any nontrivial composition of automorphisms 
$\eta(\xi^{(i)}_s)$ is distinct
from the identity, and we obtain a subgroup $\ZZ^3$ of $\pi_1(\Xi(T^5),\xi)$.
\end{proof}

The construction of Lutz \cite{L} was generalized by the author \cite{Btori} for 
$T^2$-invariant contact structures on higher-dimensional manifolds including the tori 
$T^{2n+1}$. The computation scheme of contact homology for $T^{2n+1}$ is analogous to the
case of $T^5$. However, it is necessary to study contractible closed Reeb orbits for
a $T^2$-invariant contact form on $T^{2n+1}$, in order to make sure that the cylindrical 
contact homology is well-defined. At this point, it is reasonable to formulate the 
following conjecture.

\begin{conjecture}
There exists a $T^2$-invariant contact structure $\xi$ on $T^{2n+1}$ such that the 
fundamental group $\pi_1(\Xi(T^{2n+1}),\xi)$, based at the contact structure $\xi$,
contains the subgroup $\ZZ^{2n-1}$.
\end{conjecture}

\section{Homotopy groups of higher orders}

In this section, we explain how to extend the results of section \ref{sec:pi1} to higher 
homotopy groups. 

Let $\xi_s$, $s \in [0, 1]^k$, with $\xi_s = \xi_0$ for $s \in \partial [0,1]^k$, be a 
representative of an element in $\pi_k(\Xi(M),\xi_0)$. 
We view this representative as a $(k-1)$-parameter family of loops $\xi_{s'',s'}$, 
$s'' \in \mathbb{R}/\mathbb{Z}$, for $s' \in [0,1]^{k-1}$. To each loop, we can
associate, as in section \ref{sec:pi1}, a symplectic cobordism. We therefore obtain
a $(k-1)$-parameter family of symplectic cobordisms
$(\RR \times M, \omega_{s'})$, with $s' \in [0,1]^{k-1}$, such that $\omega_{s'} = 
d(e^t \alpha_0)$ for $s' \in \partial [0, 1]^{k-1}$. 

Let $J_{s'}$ be a compatible almost complex structure on $(\RR \times M, \omega_{s'})$,
$s' \in [0,1]^{k-1}$, with the same properties as in section \ref{sec:pi1}.

This family of symplectic cobordisms may contain nongeneric $J_{s'}$-holomorphic cylinders.
Let $K_{k-1} : C_* \to C_{*+k-1}$ be the linear morphism of degree $k-1$ obtained by 
counting $J_{s'}$-holomorphic cylinders of index $1 - k$ in $(\RR \times M, \omega_{s'})$,
for $s' \in [0, 1]^{k-1}$.

\begin{proposition}
For $k > 1$, the above construction induces a group morphism
$$
\eta_k : \pi_k(\Xi(M), \xi_0) \to {\rm Mor}_{k-1}(HC_*(M,\xi_0)) .
$$
\end{proposition}

\begin{proof}
Let $K_{k-2}(i) : C_* \to C_{*+k-2}$, $i = 0, 1$, be the linear morphism obtained by 
counting $J_{s'}$-holomorphic cylinders of index $2-k$ in $(\RR \times M, \omega_{s'})$,
for $s' \in \{ i \} \times [0, 1]^{k-2}$. Counting elements in the boundary of the 
one-dimensional moduli space of $J_{s'}$-holomorphic cylinders of index $2-k$ in
$(\RR \times M, \omega_{s'})$, for $s' \in [0, 1]^{k-1}$,we obtain a relation 
$K_{k-2}(1) - K_{k-2}(0) = K_{k-1} \circ d - (-1)^{k-1} d \circ K_{k-1}$. But
$K_{k-2}(0) = K_{k-2}(1)$, since these morphisms correspond to the same $(k-2)$-parameter
family of symplectic cobordisms. Therefore, the above relation shows that $K_{k-1}$ is
a chain map.

We denote the contact homology morphism induced by $K_{k-1}$ by $\eta_k(\xi_s)$.
We need to show that this morphism is invariant under deformation of $\xi_s$.
Given a $(k+1)$-parameter family $\xi'_{s,t}$, for $s \in [0,1]^k$, $t \in [0,1]$ and satisfying 
$\xi'_{s,t} = \xi_0$ for $s \in \partial [0,1]^k$, $t \in  [0,1]$, we can construct a $k$-parameter 
family of symplectic cobordisms $(\RR \times M, \omega_{s',t})$ equipped with compatible 
almost complex structures $J_{s',t}$. The chain maps $K_k(i) : (C_*,d) \to (C_*,d)$, 
$i = 0,1$ corresponding to the symplectic cobordisms at $[0,1]^{k-1} \times \{ 0 \}$ and 
$[0,1]^{k-1} \times \{ 1 \}$ satisfy the relation $K_{k-1}(1) - K_{k-1}(0) = K_k \circ d - (-1)^k 
d \circ K_k$, where $K_k : C_* \to C_{*+k}$ is a linear morphism of degree $k$ counting 
$J_{s',t}$-holomorphic cylinders of index $-k$ in $(\RR \times M, \omega_{s',t})$, 
for $s' \in [0,1]^{k-1}$, $t \in [0,1]$. 
The relation is obtained by counting elements in the boundary of the 
one-dimensional moduli space of $J_{s',t}$-holomorphic cylinders of index $1-k$ in
$(\RR \times M, \omega_{s',t})$, for $s' \in [0, 1]^{k-1}$, $t \in [0,1]$.
It follows that the morphisms induced by $K_{k-1}(0)$ and $K_{k-1}(1)$ on contact homology 
are identical.

Finally, we show that the map $\eta$ is a group morphism. Given representatives $\xi_s$
and $\xi'_s$, $s \in [0,1]^k$ for elements in $\pi_k(\Xi(M),\xi_0)$, we represent their product
in $\pi_k(\Xi(M),\xi_0)$ by ${\bar \xi}_s$, $s = (\hat{s},\tilde{s}) \in [0,1]^{k-1} \times [0,2]$, 
where ${\bar \xi}_s = \xi_s$ if $s \in [0,1]^k$ and 
${\bar \xi}_s = \xi'_{(\hat{s},\tilde{s}-1)}$ if $s \in [0,1]^{k-1} \times [1,2]$. 
Then, the $(k-1)$-parameter family of symplectic cobordisms corresponding to the composition
$\xi_s * \xi'_s$ of $\xi_s$ and $\xi'_s$ is the union of the $(k-1)$-parameter families of symplectic 
cobordisms corresponding to $\xi_s$ and  $\xi'_s$. It then follows that 
$\eta_k(\xi_s * \xi'_s) = \eta_k(\xi_s) + \eta_k(\xi'_s)$.
\end{proof}

When $\xi_s$ is homotopic to $(\phi_s)_* \xi_0$, where $\phi_s$, $s \in [0,1]^k$, are
diffeomorphisms of $M$ such that $\phi_s = id$ for $s \in \partial [0,1]^k$,  the morphism
$\eta_k(\xi_s)$ vanishes. Indeed, if $J$ is a generic complex structure in the sense of
Lemma \ref{lem:transversal}, then the complex structures $\phi_{s*} J$ are generic as well.
Therefore, there are no holomorphic cylinders with a negative index.
In order to obtain a nonzero value for $\eta(\xi_s)$ in this case, we have to turn to an 
enhanced version of contact homology, involving $J$-holomorphic cylinders with marked points.

We obtain rigid $J$-holomorphic cylinders if we require that marked points are mapped to 
some generic cycles in the contact manifold $M$. To each such cycle in $M$ corresponds
a generator $t$ in the chain complex $C^{\bar a}_*$. 
This will be illustrated by the following result.

\begin{proposition} \label{prop:STT}
For all $n \ge 1$, the homotopy group $\pi_{2n-1}(\Xi(T^{2n} \times S^{2n-1}),\xi_0)$, 
based at the natural contact structure $\xi_0$ on the unit cotangent bundle of $T^{2n}$, 
contains an infinite cyclic subgroup.
\end{proposition}

\begin{proof}
The case $n=1$ was proved in Proposition \ref{prop:T3}, so we can assume $n > 1$.
We write $M = T^{2n} \times S^{2n-1}$. There is a natural contact form $\alpha_0$ for the contact
structure $\xi_0$ on $M$, given by
$$
\alpha_0 = \sum_{i=1}^{2n} y_i dx_i ,
$$
where $x_1, \ldots, x_{2n}$ are the coordinates on $T^{2n}$ and $\sum_{i=1}^{2n} y_i^2 = 1$.

Let $\Psi : SO(2n) \times M \to M$ be the diffeomorphism induced by  the action of $SO(2n)$ 
on the second factor of $T^{2n} \times S^{2n-1}$. Since $\pi_{2n-2}(SO(2n-1))$ is torsion
\cite{Ke}, it follows that the group morphism
$$
\pi_{2n-1}(SO(2n)) \to \pi_{2n-1}(S^{2n-1})
$$ 
in the homotopy long exact sequence induced by the fibration 
$$
SO(2n-1) \hookrightarrow SO(2n) \to S^{2n-1}
$$ 
does not vanish. Let $b : S^{2n-1} \to SO(2n)$ be a map such that the image of 
its homotopy class by the above morphism is a map $S^{2n-1} \to S^{2n-1}$ of degree
$d > 0$. Let $\xi_s = \Psi^*_{b(s)} \xi_0$, for all $s \in S^{2n-1}$.  

Consider the homotopy class ${\bar a}$ of loops wrapping once around $T^{2n}$ along the 
$x_1$-direction. There is a $(2n-1)$-torus of closed Reeb orbits in the homotopy class ${\bar a}$.
Let $D$ be the $(2n-1)$-cycle in $M$ spanning the variables $x_2, \ldots, x_{2n}$ along $T^{2n}$, 
at some value of $x_1$ and some generic point of $S^{2n-1}$. Let $t$ be the generator 
corresponding to the cycle $D$.

Let $HC^{\bar a}_*(M,\xi_0)$ be the contact homology over the coefficient ring $\QQ$, with the
additional generator $t$. Since in this case contact homology coincides with the Morse theory
of the $(2n-1)$-torus of closed Reeb orbits \cite{Bproceed}, all $J$-holomorphic cylinders are 
contained in a small neighborhood of this $(2n-1)$-torus and are disjoint from the cycle $D$. 
Therefore, $HC^{\bar a}_*(M,\xi_0)$ is isomorphic to $H_*(T^{2n-1},\QQ[t])$, with some degree shift.

With our choice of the map $b$, there are $d$ points $s_1, \ldots, s_d \in S^{2n-1}$ so that 
each closed orbit of $(\Psi_{b(s_i)})^{-1}_* R_{\alpha_0}$ in the homotopy class ${\bar a}$ intersects 
the cycle $D$ in one point, for $i = 1, \ldots, d$. In view of Lemma \ref{lem:holcyl}, there is a 
unique $J_{s_i}$-holomorphic cylinder of index $0$ passing through $D$. Adding a marked point at 
its intersection with $D$, and taking into account the condition of passing through $D$, 
we obtain a $J_{s_i}$-holomorphic cylinder of index $2-2n$.

Hence, we deduce that $\eta_{2n-1}(\xi_s)$ is the multiplication by $d \cdot t$. If we consider the
$(2n-1)$-sphere of contact structures $\xi_{g_m(s)}$, for $s \in S^{2n-1}$, where 
$g_m : S^{2n-1} \to S^{2n-1}$ is a map of degree $m$, then $\eta_{2n-1}(\xi_{g_m(s)})$ is 
the multiplication by $md \cdot t$. therefore, the $(2n-1)$-sphere $\xi_s$, $s \in S^{2n-1}$, 
has infinite order in $\pi_{2n-1}(\Xi(T^{2n} \times S^{2n-1}),\xi_0)$.
\end{proof}

Using the same ideas, we can also prove the following variant of Proposition
\ref{prop:STT}.

\begin{proposition}
Let $M$ be a compact orientable Riemannian manifold of dimension $2n$, $n \ge 1$, with 
negative sectional curvature and such that $H_1(M,\ZZ) \neq 0$. Then the homotopy group
$\pi_{2n-1}(\Xi(ST^*M),\xi_0)$, based at the natural contact structure $\xi_0$ on the unit 
cotangent bundle of $M$, contains an infinite cyclic subgroup.

\end{proposition}

\begin{proof}
The case $n=1$ was proved in Proposition \ref{prop:STSigma}, so we can assume $n>1$.
We consider a free homotopy class ${\bar a}$ in $ST^*M$ induced by a nontrivial free
homotopy class ${\bar a}_0$ in $M$. Let $\alpha_0$ be the natural contact form for $\xi_0$. 
By the assumptions on the curvature of $M$, there is a unique closed Reeb orbit in class ${\bar a}$. 
Hence, the contact homology $HC^{\bar a}_*(ST^*M,\xi_0)$ with coefficients in $\QQ$ and
any number of additional generators, is always generated by this unique closed orbit.

Since ${\bar a}_0$ induces a nonzero element of $H_1(M,\ZZ)$, we can find a cycle 
$D_0 \in H_{2n-1}(M,\ZZ)$ such that the intersection number of ${\bar a}_0$ and $D_0$ 
is equal to $d \neq 0$. Let $D \in H_{2n-1}(ST^*M;\ZZ)$ be a lift of the cycle $D_0$, so that the
preimage in $ST^*M$ of the closed geodesic in class ${\bar a}_0$ has intersection number $d$
with $D$.

The definition of $\xi_s$, $s \in S^{2n-1}$, and the remainder of the proof are the same as in 
Proposition \ref{prop:STT}.
\end{proof}

\section*{Appendix : proof of Lemma \ref{lem:transversal}}

Lemma \ref{lem:transversal} was proved independently by Dragnev \cite{D} using a slightly 
different argument.
We will follow the arguments from \cite{MS}, and adapt them to this situation.
Let $\mathcal{J}^l(M,\alpha)$ be the set of $\RR$-invariant compatible almost complex 
structures of class $C^l$ on $(\RR \times M, d(e^t\alpha))$ satisfying $J\xi = \xi$ and 
$J\frac{\partial}{\partial t} = R_\alpha$. It is a smooth Banach manifold, and its tangent space
$T_J\mathcal{J}^l(M,\alpha)$ is the Banach space of endomorphisms $Y$ of $\xi$ over $M$, of 
class $C^l$, such that $YJ + JY = 0$ and $d\alpha(\cdot,Y\cdot) + d\alpha(Y\cdot,\cdot)= 0$.

Let $\og$ and $\ug$ be closed Reeb orbits with periods $\overline{T}$ and $\underline{T}$
in class ${\bar a}$, let $d >0$, $1 \le k \le l$ and $p > 2$. Let us introduce coordinates \
$(\vartheta,z) \in
\RR/\ZZ \times D^{2n-2}$ in a tubular neighborhood of these closed Reeb orbits, so that
$\RR/\ZZ \times \{ 0 \}$ parametrizes the closed obits. We define $\mathcal{B}^{p,d}_k(\og,\ug)$ as 
the set of maps $F = (a,f) : \RR \times S^1 \to \RR \times M$, locally in $L^p_k$, such that
\begin{eqnarray*}
\lim_{s \to +\infty} a(s,\theta) &=& +\infty , \\
\lim_{s \to +\infty} f(s,\theta) &=& \og(\overline{T}(\theta-\overline{\theta}_0)), 
\qquad \textrm{for some } \overline{\theta}_0 \in S^1, \\
\lim_{s \to -\infty} a(s,\theta) &=& -\infty, \\
\lim_{s \to -\infty} f(s,\theta) &=& \ug(\underline{T}(\theta-\underline{\theta}_0)) ,
\qquad \textrm{for some } \underline{\theta}_0 \in S^1, 
\end{eqnarray*}
the maps $a(s,\theta) - \overline{T}s - \overline{a}_0$, $\vartheta(s,\theta) - 
(\theta-\overline{\theta}_0)$
and $z(s,\theta)$ are in 
$$
L_k^{p,d}(\RR^+ \times S^1) = L_k^p(\RR^+ \times S^1, e^{ds} ds \, d\theta) 
$$
for some $\overline{a}_0 \in \RR$, and the maps $a(s,\theta) - \underline{T}s - \underline{a}_0$, 
$\vartheta(s,\theta) - (\theta-\underline{\theta}_0)$ and $z(s,\theta)$ are in 
$$
L_k^{p,d}(\RR^- \times S^1) = L_k^p(\RR^+ \times S^1, e^{-ds} ds \, d\theta) 
$$
for some $\underline{a}_0 \in \RR$.

Let $\beta : \RR \to \RR$ be a smooth function such that $\beta(s) = 0$ if $s \le 0$, $\beta(s) = 1$
if $s \ge 1$, and $0 \le \beta'(s) \le 2$. For $F \in \mathcal{B}^{p,d}_k(\og,\ug)$, let $\overline{V}$ 
be the vector space generated by the sections $\beta(s) R_\alpha$ and 
$\beta(s) \frac{\partial}{\partial t}$ of $F^*T(\RR \times M)$; let $\underline{V}$ be the vector space 
generated by the sections $\beta(-s) R_\alpha$ and $\beta(-s) \frac{\partial}{\partial t}$ 
of $F^*T(\RR \times M)$.

With these definitions, $\mathcal{B}^{p,d}_k(\og,\ug)$ is a smooth Banach manifold and its 
tangent space $T_F \mathcal{B}^{p,d}_k(\og,\ug)$ is the Banach space
$\overline{V} \oplus L_k^{p,d}(\RR \times S^1; F^*T(\RR \times M)) \oplus \underline{V}$.

Let $\mathcal{M}^l(\og,\ug)$ be the universal moduli space, consisting of pairs $(F,J) \in
\mathcal{B}^{p,d}_k(\og,\ug) \times \mathcal{J}^l(M,\alpha)$ such that $F$ is a $J$-holomorphic
cylinder, for the complex structure $j$ on $\RR \times S^1$ such that $j \frac{\partial}{\partial s}
= \frac{\partial}{\partial \theta}$. Since $J$-holomorphic cylinders converge exponentially
to closed Reeb orbits \cite{HWZ}, all holomorphic cylinders are in $\mathcal{B}^{p,d}_k(\og,\ug)$
when $d > 0$ is sufficiently small.

For $J \in \mathcal{J}^l(M,\alpha)$, let $g_J$ be the $J$-invariant Riemannian metric on 
$\RR \times M$ defined by
$$
g(v,w) = d\alpha(v,Jw) + \alpha(v) \alpha(w) + dt(v) dt(w) .
$$
Linearizing the Cauchy-Riemann operator at $(F,J)$ using the exponential map of $g_J$, 
we obtain 
\begin{eqnarray*}
L_{(F,J)}\delbar :  T_F \mathcal{B}^{p,d}_k(\og,\ug) \oplus T_J\mathcal{J}^l(M,\alpha) &\to&
L_{k-1}^{p,d}(\RR \times S^1;\Lambda^{0,1}(F^*T(\RR \times M)))  \\
(\zeta,Y) &\mapsto& L_F\delbar_J \zeta + Y \circ df \circ j ,
\end{eqnarray*}
The linear map $L_F\delbar_J$ is Fredholm and has the form
$$
(L_F\delbar_J \zeta) (\frac{\partial}{\partial s}) = \nabla_\frac{\partial F}{\partial s} \zeta
+ J \nabla_\frac{\partial F}{\partial \theta} \zeta + (\nabla_\zeta J) \frac{\partial F}{\partial \theta} ,
$$
where $\nabla$ denotes the Levi-Civita connection of $g_J$.

We want to show that $\mathcal{M}^l(\og,\ug)$ is a Banach manifold. Using an infinite dimensional
implicit function theorem, it suffices to show that the linearized Cauchy-Riemann operator 
$L_{(F,J)}\delbar$ is surjective for all $(F,J) \in \mathcal{M}^l(\og,\ug)$. Then, the generic complex
structures can be obtained by the Sard-Smale theorem as the regular values of the natural 
projection map $\pi : \mathcal{M}^l(\og,\ug) \to \mathcal{J}^l(M,\alpha)$.
We prove surjectivity of $L_{(F,J)}\delbar$ in the case $k=1$. The case $k > 1$ can then be
deduced by elliptic regularity as in \cite{MS}.

First assume that $F$ is a vertical cylinder over a closed Reeb orbit $\gamma = 
\og = \ug$. In this case, the operator $L_{(F,J)}\delbar$ splits in two summands. The first summand
is an $\RR$-invariant operator acting on sections of $\xi$, with nondegenerate asymptotics. It is 
well-known (see for example \cite{S}) that such an operator is an isomorphism. 
The second summand is the actual Cauchy-Riemann operator acting on sections of 
$\CC = \RR \{\frac{\partial}{\partial t},R_\alpha\}$, over $\CC^*$. Its index is $2$ and its 
kernel consists of constant sections. Hence, both summands are surjective, and so 
is $L_{(F,J)}\delbar$ .

Let us argue by contradiction for the general case $\og \neq \ug$. By the Hahn-Banach theorem,
we can find $\eta \neq 0$ in $L^{q,d}(\RR \times S^1;\Lambda^{1,0}(F^*T(\RR \times M)))$, 
with $\frac1p + \frac1q = 1$, such that
$$
\langle \eta, L_F\delbar_J \zeta \rangle + \langle \eta, Y \circ df \circ j \rangle = 0
$$
for all $(\zeta,Y) \in T_F \mathcal{B}^{p,d}_1(\og,\ug) \oplus T_J\mathcal{J}^l(M,\alpha)$.
Since $\zeta$ and $Y$ can be chosen independently, the above two terms have to vanish
separately.

Let $p \in \RR \times S^1$ be a point such that
$$
df_p \neq 0 \quad \textrm{and} \quad f^{-1}(f(p)) = \{ p \} .
$$
Such a point $p$ is called an injective point.
Choose $Y \in T_J\mathcal{J}^l(M,\alpha)$ with support in a small ball centered at $F(p)$. 
Since $df \circ j$ is injective at $p$ and $Y_{F(p)}$ can be any endomorphism of $\xi_{F(p)}$, 
it follows from
$$
\langle \eta, Y \circ df \circ j \rangle = 0
$$
that the orthogonal projection of $\eta$ to $\xi$ vanishes at $p$ : $\eta_\xi(p) = 0$. 

Let us prove that, for $R > 0$ large enough, the set of injective points $p \in \RR \times S^1$  
is open and dense in $(R, +\infty) \times S^1$. It will follow that $\eta_\xi = 0$ 
on $(R,+\infty) \times S^1$. 

\begin{lemma}  \label{lem:embed}
If $R > 0$ is sufficiently large, then the map $f$ restricted to $(R,\infty) \times S^1$
is an embedding.
\end{lemma}

\begin{proof}
By Theorem 2.8 of \cite{HWZ}, there are local coordinates
$(\vartheta,z) \in S^1 \times \RR^{2n-2}$ parametrizing a tubular neighborhood of $\og$ so
that 
\begin{equation} \label{z}
z(s,\theta) = e^{\int_{s_0}^s \gamma(\tau) d\tau} [e(\theta) + r(s,\theta)]
\end{equation}
with $\gamma(s) \to \lambda < 0$ as $s \to +\infty$, $e(\theta) \neq 0$ for all $\theta \in S^1$, and
all derivatives of $r(s,\theta)$ uniformly converging to $0$ as $s \to +\infty$.
Moreover, the component $\vartheta(s,\theta)$ and its derivatives satisfy
$$
| \partial^I [\vartheta(s,\theta) - \theta] | \le C_I e^{-ds}
$$
for any multi-index $I$ and some $C_I, d > 0$.

In particular, if $s$ is sufficiently large, we have
\begin{eqnarray}
\| \frac{\partial \vartheta}{\partial \theta}(s,\theta) \| &\ge& \frac12, \label{thetat}\\
\| \frac{\partial z}{\partial s}(s,\theta)\| &\ge& C e^{\int_{s_0}^s\gamma(\tau)d\tau}, 
\label{zs} \\ 
\| \frac{\partial^{i_1 + i_2} \vartheta}{\partial s^{i_1} \partial \theta^{i_2}}(s,\theta) \| &\le& 
C e^{-ds}, \label{theta2} \\
\| \frac{\partial^{i_1 + i_2} z}{\partial s^{i_1} \partial \theta^{i_2}}(s,\theta) \| &\le& 
C e^{\int_{s_0}^s\gamma(\tau)d\tau} \label{z2} 
\end{eqnarray}
with $i_1, i_2 \ge 0$ and $i_1+i_2 \ge 2$.
It follows from (\ref{thetat}) and (\ref{zs}) that, if $R$ is sufficiently large, 
the map $f|_{(R,\infty) \times S^1}$ is an immersion.

Let us show that, for $R$ sufficiently large, the map $f|_{(R,\infty)\times S^1}$ is injective.
Assume by contradiction that there exist sequences $(s_n,\theta_n)$ and $(s'_n,\theta'_n)$ 
with $s_n, s'_n \to \infty$, $\theta_n, \theta'_n \in S^1$ and $f(s_n,\theta_n) = f(s'_n,\theta'_n)$.
We write $(s'_n-s_n, \theta'_n-\theta_n) = \rho_n {\vec v}_n$ with $\rho_n \in \RR^+$ and 
$\| {\vec v}_n \| = 1$.

The second order Taylor expansion of $\vartheta$ along the line joining
$(s_n,\theta_n)$ and $(s'_n,\theta'_n)$ reads
$$
\vartheta(s'_n,\theta'_n) = \vartheta(s_n,\theta_n) 
+ \rho_n \frac{\partial \vartheta}{\partial {\vec v}_n}(s_n,\theta_n) 
+ \frac{\rho_n^2}2 \frac{\partial^2 \vartheta}{\partial {\vec v}_n^2}(s''_n,\theta''_n)
$$
for some $(s''_n,\theta''_n)$ between $(s_n,\theta_n)$ and $(s'_n,\theta'_n)$, so that
$$
\frac{\partial \vartheta}{\partial {\vec v}_n}(s_n,\theta_n) 
+ \frac{\rho_n}2 \frac{\partial^2 \vartheta}{\partial {\vec v}_n^2}(s''_n,\theta''_n) =0.
$$
In view of (\ref{theta2}), the second term converges to zero so that ${\vec v}_n \to (\pm 1, 0)$.
In particular, $t'_n - t_n \to 0$ as $n \to \infty$.

The second order Taylor expansion of a component $z_i$ of $z$ along the same line gives similarly
$$
\frac{\partial z_i}{\partial {\vec v}_n}(s_n,\theta_n) 
+ \frac{\rho_n}2 \frac{\partial^2 z_i}{\partial {\vec v}_n^2}(s^{(i)}_n,\theta^{(i)}_n) = 0
$$
for some $(s^{(i)}_n,\theta^{(i)}_n)$ between $(s_n,\theta_n)$ and $(s'_n,\theta'_n)$.
Therefore, using (\ref{zs}) and (\ref{z2}), it follows that $\rho_n$ is bounded away from zero. 

On the other hand, since $z(s_n,\theta_n) = z(s'_n,\theta'_n)$ it follows from (\ref{z}) that
$$
e^{\int_{s'_n}^{s_n} \gamma(\tau) d\tau} [e(\theta_n) +r(s_n,\theta_n)] = e(\theta_n) 
+ [e(\theta'_n)-e(\theta_n)] + r(s'_n,\theta'_n).
$$
In this expression, $r(s_n,\theta_n)$, $e(\theta'_n)-e(\theta_n)$ and $r(s'_n,\theta'_n) \to 0$ 
as $n \to \infty$, while $e(\theta_n)$ is bounded away from zero. 
Hence, $s'_n - s_n \to 0$ as well, contradicting the fact that $\rho_n$ is bounded away from zero.
\end{proof}

Since $\og \neq \ug$, it is clear that 
$f((R,+\infty) \times S^1) \cap f ((-\infty,-R) \times S^1) = \emptyset$ for $R > 0$
sufficiently large.
In view of Lemma \ref{lem:embed}, we just have to consider the intersection 
$I = f([-R,R] \times S^1) \cap f((R,+\infty) \times S^1)$. Since $f([-R,R] \times S^1)$  is compact, 
$f^{-1}(I)$ is closed in $(R,+\infty) \times S^1$.  Suppose $f^{-1}(I) \cap (R,+\infty) \times S^1$ 
contains an open subset $U_1$, and let $U_2 = f^{-1}(f(U_1)) \cap (-R,R) \times S^1$.
Since $f |_{U_1}$ is bijective, we can define $\phi = (f |_{U_1})^{-1} \circ f|_{U_2} : U_2 \to U_1$.
With this map, we have $f(p) = f(\phi(p))$ for all $p \in U_2$. Now, $f^* d\alpha$ cannot vanish on 
an open set : by the unique continuation theorem $F$ would then be a cylinder over a Reeb orbit.
Therefore, the map $\phi$ preserves the conformal structure of $\RR \times S^1$ on a dense
subset, and hence on the cylinder. Since the function $a$ is determined 
by the Cauchy-Riemann equation up to an additive constant, then for some constant $C \in \RR$,
we have $a(p) = a(\phi(p)) + C$, for all $p \in U_2$. 
Therefore, we can apply the unique continuation theorem to $F$ and $F \circ \phi$. 
We obtain a sequence $p_i = (s_i,\theta_i) \in \RR \times S^1$ satisfying $s_i \to +\infty$ and
$f(p_i) = f(p)$ for some $p \in U_2$. This contradicts the (exponential) convergence of $f$ to $\og$
as $s \to +\infty$. Therefore, the closed set $f^{-1}(I) \cap (R,+\infty) \times S^1$ cannot 
contain an open set, and the set of injective points in open and dense in $\RR \times S^1$.

We now want to show that $\eta_\xi = 0$ in $(R,+\infty) \times S^1$ implies that $\eta = 0$ on 
an open subset of $(R,+\infty) \times S^1$. Let $v$ be the orthogonal projection of 
$\frac{\partial F}{\partial s}$ (or of $\frac{\partial f}{\partial s}$) to $\xi$. Then, the orthogonal 
projection of $\frac{\partial F}{\partial \theta}$ to $\xi$ is $Jv$. Let us compute 
$g_J(R_\alpha, (L_F\delbar_J \zeta) (\frac{\partial}{\partial s}))$, when $\zeta$ is a section of
$F^* \xi$ :
\begin{eqnarray*}
g_J(R_\alpha, (L_F\delbar_J \zeta) (\frac{\partial}{\partial s}))
&\!=&\! g_J(R_\alpha, \nabla_\frac{\partial F}{\partial s} \zeta) 
+ g_J(R_\alpha, J \nabla_\frac{\partial F}{\partial \theta} \zeta) \\
&& + g_J(R_\alpha,(\nabla_\zeta J)\frac{\partial F}{\partial \theta}) \\
&\!=&\! g_J(R_\alpha, \nabla_v \zeta) 
- g_J(J R_\alpha, \nabla_\frac{\partial F}{\partial \theta} \zeta) 
+ g_J(R_\alpha,(\nabla_\zeta J) Jv) \\
&\!=&\! g_J(R_\alpha, \nabla_v \zeta) + g_J(R_\alpha,\nabla_\zeta ( J^2 v))
- g_J(R_\alpha, J \nabla_\zeta (Jv)) \\
&\!=&\! g_J(R_\alpha, \nabla_v \zeta) - g_J(R_\alpha, \nabla_\zeta v)
+ g_J(J R_\alpha, \nabla_\zeta (Jv)) \\
&\!=&\! g_J(R_\alpha, [v, \zeta])  \\
&\!=&\! \alpha([v,\zeta]) \\
&\!=&\! -d\alpha(v,\zeta) .
\end{eqnarray*}
We obtain
$$
g_J(\frac{\partial}{\partial t}, (L_F\delbar_J \zeta) (\frac{\partial}{\partial s})) =
d\alpha(Jv,\zeta)
$$
in a similar way. We already saw that $v$ cannot vanish on an open set, so we can restrict
ourselves to points where $v \neq 0$. Since $v$ and $Jv$ are then linearly independent and 
$d\alpha$ is nondegenerate on $\xi$,  for $p \in (R,+\infty) \times S^1$ and $\eta(p) \neq 0$ we can
choose $\zeta(p)$ in such a way that the pairing between $\eta(p)$ and $L_F\delbar_J \zeta(p)$ 
does not vanish. Taking $\zeta$ to be a section of $F^*\xi$ with small support around $p$,
the equation
$$
\langle \eta, L_F\delbar_J \zeta \rangle = 0
$$
gives $\eta(p) = 0$.  Since $p$ is arbitrary in $(R,+\infty) \times S^1$, $\eta$ vanishes in a 
neighborhood of $p$.

Since $\eta$ is in the kernel of the formal adjoint $(L_{(F,J)}\delbar)^*$, we can apply the unique
continuation theorem, so that $\eta = 0$ on $\RR \times S^1$.
This contradicts our assumption that $\eta \neq 0$ annihilates the image of the linearized operator 
$L_{(F,J)}\delbar$, and completes the proof of Lemma \ref{lem:transversal}.


\begin{thebibliography}{99}
\bibitem{Btori} F. Bourgeois, Odd dimensional tori are contact manifolds, 
{\it Int. Math. Res. Not.} {\bf 2002} No. 30, 1571--1574.

\bibitem{Bproceed} F. Bourgeois, A Morse-Bott approach to contact homology, 
in ``Symplectic and Contact Topology : Interactions and Perspectives'', {\it
Fields Institute Communications} {\bf 35} (2003), 55--77.

\bibitem{BM} F. Bourgeois, K. Mohnke, Coherent orientations in Symplectic Field Theory,
{\it Math. Z.} {\bf 248} (2004), No. 1, 123--146.

\bibitem{D} D. Dragnev, Fredholm theory and transversality for noncompact pseudoholomorphic 
maps in symplectizations,  {\it Comm. Pure Appl. Math.}  {\bf 57}  (2004) No. 6, 726--763.

\bibitem{EGH} Y.~Eliashberg, A.~Givental and H.~Hofer, Introduction to Symplectic Field 
Theory, {\it Geom. Funct. Anal., Special Volume, Part II} (2000), 560--673.

\bibitem{FH} A. Floer, H. Hofer, Coherent orientations for periodic orbit problems in symplectic
geometry, {\it Math. Z.} {\bf 212} (1993) No. 1, 13--38.

\bibitem{GG} H. Geiges, J. Gonzalo, On the topology of the space of contact structures 
on torus bundles, {\it Bull. London Math. Soc.} {\bf 36} (2004), 640--646.


\bibitem{G} E. Giroux, Une infinit\'e de structures de contact tendues sur une infinit\'e
de vari\'et\'es, {\it Invent. Math.} {\bf 135} (1999), 789--802.

\bibitem{HWZ} H. Hofer, K. Wysocki, E. Zehnder, Pseudoholomorphic curves in 
symplectizations I : Asymptotics, {\it Ann.~Inst.~Henri Poincar\'{e}, Analyse Nonlin\'{e}aire}
{\bf 13} (1996) No. 3, 337--379.

\bibitem{K} T. K\'alm\'an, Contact homology and one parameter families of Legendrian
knots, {\it arXiv preprint} (math.GT/0407347).

\bibitem{Kanda} Y. Kanda, The classification of tight contact structures on the $3$-torus,
{\it Comm. Anal. Geom.} {\bf 5} (1997) No. 3, 413--438.

\bibitem{Ke} M.A. Kervaire, Some nonstable homotopy groups of Lie groups, {\it Illinois J. Math.}
{\bf 4} (1960), 161--169.

\bibitem{L} R. Lutz, Sur la g\'eom\'etrie des structures de contact invariantes, 
{\it Ann. Inst. Fourier (Grenoble)} {\bf 29} (1979), 283--306.

\bibitem{MS} D. McDuff, D. Salamon, $J$-holomorphic curves and Quantum Cohomology, 
{\em University Lecture Series} {\bf 6}, AMS, 1994.

\bibitem{S} D. Salamon, Lectures on Floer homology, in ``Symplectic Geometry and Topology'',
{\it IAS/Park City Mathematics Series} {\bf 7} (1999), 145--229.
\end{thebibliography}
\end{document}